\documentclass[reqno]{amsart}
\usepackage{amsmath}
\allowdisplaybreaks[3]
\numberwithin{equation}{section}
\numberwithin{equation}{section}

\def\proof{\indent{\em Proof.\quad}}
\def\endproof{\hfill\hbox{$\sqcup$}\llap{\hbox{$\sqcap$}}\medskip}
\newtheorem{thm}{\indent Theorem}[section]
\newtheorem{cor}[thm]{\indent Corollary}
\newtheorem{lem}[thm]{\indent Lemma}
\newtheorem{prop}[thm]{\indent Proposition}
\newtheorem{dfn}{{\indent\bf Definition}}[section]
\newtheorem{rmk}{{\indent\bf Remark}}[section]
\newtheorem{expl}{{\indent\bf Example}}[section]

\newcommand{\mb}{\mbox}
\newcommand{\hs}{\hspace}

\newcommand{\ol}{\overline}

\newcommand{\strl}{\stackrel}

\newcommand{\td}{\tilde}

\newcommand{\fr}{\frac}
\newcommand{\ed}{{\rm End}}
\newcommand{\edd}{\end{document}}
\newcommand{\be}{\begin{equation}}
\newcommand{\ee}{\end{equation}}

\newcommand{\lmx}{\left(\begin{matrix}}
\newcommand{\rmx}{\end{matrix}\right)}
\newcommand{\ldt}{\left|\begin{matrix}}
\newcommand{\rdt}{\end{matrix}\right|}

\newcommand{\tr}{{\rm tr\,}}
\newcommand{\vfi}{\varphi}

\newcommand{\bbr}{{\mathbb R}}

\newcommand{\bbc}{{\mathbb C}}

\newcommand{\ba}{\begin{array}}
\newcommand{\ea}{\end{array}}
\newcommand{\nnm}{\nonumber}
\newcommand{\beal}{\begin{align}}
\newcommand{\eal}{\end{align}}
\newcommand{\bea}{\begin{eqnarray}}
\newcommand{\eea}{\end{eqnarray}}

\newcommand{\ppp}[3]{\fr{\partial^2 #1}{\partial #2\partial #3}}
\newcommand{\dd}[2]{\fr{d #1}{d #2}}
\begin{document}

\title[On the equiaffine symmetric hyperspheres]{On the equiaffine symmetric hyperspheres}\thanks{Research supported by
NSFC (No. 11171091, 11371018) and partially supported by NSF of Henan Province (no. 132300410141).}%

\author{Xingxiao Li\,$^*$}\thanks{$^*$ The corresponding author}
\author{Guosong Zhao}%
\date{}

\begin{abstract}
We introduce and study the equiaffine symmetric {\bf hyperspheres}. For the first step we consider the locally strongly convex ones. In fact, by the idea used by Naitoh, we provide in this paper a direct proof of the complete classification for those affine symmetric hyperspheres. Then, via an earlier result of the first author, we are able to provide an alternative proof for the classification theorem of the affine hypersurface with parallel Fubini-Pick forms, which has already been established by Z.J. Hu et al in a totally different way.

\vskip 0.1in
{\bf Key words and expressions}\ equiaffine {\bf hypersphere}, affine metric, Fubini-Pick form, symmetric space, affine symmetric hypersurface

\vskip 0.1in
{\bf 2000 AMS classification:} Primary 53A15; Secondary 53B25
\end{abstract}
\maketitle

\tableofcontents

\section{Introduction}

As we know, affine hyperspheres are very special in the equiaffine differential geometry of hypersurfaces. In particular, if an affine hypersurface is of parallel Fubini-Pick form, then it must be an affine hypersphere (\cite{bok-nom-sim90}). If we only take account of the definition, affine hyperspheres seem very simple but in fact they do form a very large class of hypersurfaces. Consequently it is a great challenge to find explicitly all the affine hyperspheres and now it still remains a very hard job. Although this, the study of affine hyperspheres has been made a lot of great achievement by many authors. For example, the proof of the Calabi's conjecture (see for example, \cite{amli90}, \cite{amli92}), the classification of hyperspheres of constant sectional affine curvatures (\cite{vra-li-sim91},  \cite{wang93} and \cite{kri-vra99}), the generalizations of Calabi's composition of affine hyperbolic hyperspheres (with multiple factors, \cite{lix93}; for more general cases, \cite{dil-vra94}), the characterization of the Calabi's composition of hyperbolic hyperspheres (\cite{hu-li-vra08}; also \cite{lix13} and \cite{lix14} in a different manner), and the classification of locally strongly convex hypersurfaces with parallel Fubini-Pick forms (\cite{dil-vra-yap94} and \cite{hu-li-sim-vra09} for some special cases; \cite{hu-li-vra11} for general case). As for the general nondegenerate case, there also have been some interesting partial classification results, see for example the series of published papers by Z.J. Hu et al: \cite{hu-li11}, \cite{hu-li-li-vra11a} and \cite{hu-li-li-vra11b}. In this direction, a very recent development is the preprint article \cite{Hil12} in which the author aimed at a complete classification of nondegenerate centroaffine hypersurfaces with parallel Fubini-Pick form.

In this paper, on the basis of a recent characterization of Calabi composition of hyperbolic hypersphere (\cite{lix13}, \cite{lix14}), we make use of the idea by H. Naitoh  in \cite{nai81} for classification of totally real parallel submanifolds in the complex projective space, to provide a direct proof of the complete classification of symmetric affine hyperspheres. Then, via an earlier result of the author, we easily give an alternative and simpler proof for the classification theorem (Theorem \ref{cla thm}) for the affine hypersurface with parallel Fubini-Pick forms, which has already been established by Z.J. Hu et al in a totally different way (see \cite{hu-li-vra11} for the detail).

Our main theorem is stated as follows:

{\thm[The main theorem]\label{main} Let $x:M^n\to \bbr^{n+1}$ ($n\geq 2$) be a locally strongly convex affine {\bf hypersphere}. If $x$ is locally affine symmetric, then either of the following two cases holds:

$(1)$ With the affine metric $g$, the Riemannian manifold $(M^n,g)$ is irreducible and $x$ is locally affine equivalent to

$(a)$ one of the three kinds of quadric affine spheres: Ellipsoid, elliptic paraboloid and hyperboloid; or

$(b)$ the standard embedding of the Riemannian symmetric space ${\rm SL}(m,\bbr)/{\rm SO}(m)$ into $\bbr^{n+1}$ with $n=\fr12m(m+1)-1$, $m\geq 3$;
or

$(c)$ the standard embedding of the Riemannian symmetric space ${\rm SL}(m,\bbc)/{\rm SU}(m)$ into $\bbr^{n+1}$ with $n=m^2-1$, $m\geq 3$; or

$(d)$ the standard embedding of the Riemannian symmetric space ${\rm SU}^*(2m)/{\rm Sp}(m)$ into $\bbr^{n+1}$ with $n=2m^2-m-1$, $m\geq 3$; or

$(e)$ the standard embedding of the Riemannian symmetric space ${\rm E}_{6(-26)}/{\rm F}_4$ into $\bbr^{27}$.

$(2)$ $(M^n,g)$ is reducible and $x$ is locally affine equivalent to the Calabi product of $r$ points and $s$ of the above irreducible hyperbolic affine spheres of lower dimensions, where $r$, $s$ are nonnegative integers and $r+s\geq 2$.}

Examples (b), (c), (d) and (e) are explicitly presented in Section 3, while examples in (a) can be found in the most text books, see for example \cite{li-sim-zhao93}.

{\sc Acknowledgement} The first author is grateful to Professor A-M Li for his encouragement and important suggestions during the preparation of this article. He also thanks Professor Z.J. Hu for providing him valuable related references some of which are listed in the end of this paper.

\section{Preliminaries}

\subsection{The equiaffine geometry of hypersurfaces}

In this subsection, we brief some basic facts in the equiaffine geometry of hypersurfaces. For details the readers are referred to some text books, say, \cite{li-sim-zhao93} and \cite{nom-sas94}.

Let $x:M^n\to\bbr^{n+1}$ be a nondegenerate hypersurface. Then there are several basic equiaffine invariants of $x$ among which are: the affine metric (Berwald-Blaschke metric) $g$, the affine normal $\xi:=\fr1n\Delta_gx$, the Fubini-Pick $3$-form (the so called cubic form) $A\in\bigodot^3T^*M^n$ and the affine second fundamental $2$-form $B\in\bigodot^2T^*M^n$. By using the index lifting by the metric $g$, we can identify $A$ and $B$ with the linear maps $A:TM^n\to \ed(TM^n)$ or $A:TM^n\bigodot TM^n\to TM^n$ and $B:TM^n\to TM^n$, respectively, by
\be\label{ab}
g(A(X)Y,Z)=A(X,Y,Z) \mb{\ or\ }g(A(X,Y),Z)=A(X,Y,Z),\quad
g(B(X),Y)=B(X,Y),
\ee
for all $X,Y,Z\in TM^n$. Sometimes we call the corresponding $B\in \ed(TM^n)$ the affine shape operator of $x$. In this sense, the affine Gauss equation can be written as follows:
\be\label{gaus}
R(X,Y)Z=\fr12(g(Y,Z)B(X)+B(Y,Z)X-g(X,Z)B(Y)-B(X,Z)Y)-[A(X),A(Y)](Z),
\ee
where, for any linear transformations $T,S\in \ed(TM^n)$,
\be\label{comm}
[T,S]=T\circ S-S\circ T.
\ee
Each of the eigenvalues $B_1,\cdots,B_n$ of the affine shape operator $B:TM^n\to TM^n$ is called the affine principal curvature of $x$. Define
\be\label{afme}
L_1:=\fr1n\tr B=\fr1n\sum_iB_i.
\ee
Then $L_1$ is referred to as the affine mean curvature of $x$. A hypersurface $x$ is called an (elliptic, parabolic, or hyperbolic) affine hypersphere, if all of its affine principal curvatures are equal to one (positive, 0, or negative) constant. In this case we have
\be\label{afsp}
B(X)=L_1X,\quad\mb{for all\ }X\in TM^n.
\ee
It follows that the affine Gauss equation \eqref{gaus} of an affine hypersphere assumes the following form:
\be\label{gaus_af sph}
R(X,Y)Z=L_1(g(Y,Z)X-g(X,Z)Y)-[A(X),A(Y)](Z),
\ee

Furthermore, all the affine lines of an elliptic affine hypersphere or a hyperbolic affine hypersphere $x:M^n\to\bbr^{n+1}$ pass through a fix point $o$ which is refer to as the affine center of $x$; Both the elliptic affine hyperspheres and the hyperbolic affine hyperspheres are called proper affine hyperspheres, while the parabolic affine hyperspheres are called improper affine hyperspheres.

For each vector field $\eta$ transversal to the tangent space of $x$, we have the following direct decomposition of vector spaces
$$
x^*T\bbr^{n+1}=x_*(TM)+\bbr\cdot\eta.
$$
This decomposition and the canonical differentiation $\bar D^0$ on $\bbr^{n+1}$ define a nondegenerate bilinear form $h\in\bigodot^2T^*M^n$ and a connection $D^\eta$ on $TM^n$ as follows:
\be\label{dfn h}
\bar D^0_XY=x_*(D^\eta_XY)+h(X,Y)\eta,\quad\forall X,Y\in TM^n.
\ee
\eqref{dfn h} can be referred as to the affine Gauss formula of the hypersurface $x$.

In what follows we make the following convention for the range of indices:
$$1\leq i,j,k,l\leq n.$$

Let $\{e_i,e_{n+1}\}$ be a local unimodular frame field along $x$, and $\{\omega^i,\omega^{n+1}\}$ its dual coframe.  Then $\eta:=e_{n+1}$ is transversal to the tangent space $x_*(TM)$. Write $h=\sum h_{ij}\omega^i\omega^j$ with $h_{ij}=h(e_i,e_j)$ and $H=|\det(h_{ij})|$. Then the locally defined nondegenerate metric $g:=H^{-\fr1{n+2}}h$ is independent of the choice of the unimodular frame field $\{e_i,e_{n+1}\}$ and thus is in fact a globally well-defined metric on $M^n$ which is called the affine (or Berwald-Blaschke) metric. By taking $x$ as an $\bbr^{n+1}$-valued smooth function on $M^n$, we call the vector function $\xi:=\fr1n\tr_g(x)$ the affine normal vector.

If, in particular, $\eta$ is chosen to be parallel to the affine normal $\xi$, Then the induced connection $\nabla:=D^\eta$ is independent of the choice of $\eta$ and is referred to as the affine connection of $x$.
If $\hat\nabla$ is the Levi-Civita connection of the affine metric $g$, then
the Fubini-Pick form (as a symmetric $(1,2)$ tensor) is defined by
\be\label{f-p}
A(X,Y)=\nabla_XY-\hat\nabla_XY,\quad \forall\, X,Y\in TM,
\ee
which is identified via the affine metric $g$ with a symmetric cubic form $A(X,Y,Z)=g(A(X,Y),Z)$. This cubic form $A$ is also referred to as the Fubini-Pick form.

From now on we assume that the transversal vector $e_{n+1}$ above is parallel to the normal vector $\xi$. Then it holds that $\xi=H^{\fr1{n+2}}e_{n+1}$ and we have connection forms $\omega^A_B$, $1\leq A,B\leq n+1$, defined by
$$
d\omega^A=\omega^B\wedge\omega^A_B,\quad d\omega^A_B=\sum_{C=1}^{n+1}\omega^C_B\wedge\omega^A_C,\quad \omega^{n+1}\equiv 0.
$$
Furthermore, the local expressions of $g$, $A$ and $B$:
\be\label{gab}
A=\sum A_{ijk}\omega^i\omega^j\omega^k,\quad B=\sum B_{ij}\omega^i\omega^j,
\ee
is subject to the following basic formulas:
\begin{align}
&\sum_{i,j} g^{ij}A_{ijk}=0\text{\ (the apolarity)},\label{basic1}\\
&A_{ijk,l}-A_{ijl,k}=\fr12(g_{ik}B_{jl}+g_{jl}B_{ik} -g_{il}B_{jk}-g_{jk}B_{il}),\label{basic3}\\
&\sum_{l}A^l_{ij,l}=\fr n2(L_1g_{ij}-B_{ij}),\label{basic3-1}
\end{align}
where $A_{ijk,l}$ are the covariant derivatives of $A_{ijk}$ with respect to the Levi-Civita connection of $g$.

Define
\be\label{hijk0}
\sum_kh_{ijk}\omega^k=dh_{ij}+h_{ij}\omega^{n+1}_{n+1}-\sum h_{kj}\omega^k_i-\sum h_{ik}\omega^k_j.
\ee
Then the Fubini-Pick form $A$ can be determined by the following formula:
\be\label{hijktoaijk}
A_{ijk}=-\fr12H^{-\fr1{n+2}}h_{ijk}.
\ee

Define the normalized scalar curvature $\chi$ and the Pick invariant $J$ by
$$
\chi=\fr1{n(n-1)}\sum g^{il}g^{jk}R_{ijkl},\quad J=\fr1{n(n-1)}\sum A_{ijk}A_{pqr}g^{ip}g^{jq}g^{kr}.$$
Then the affine Gauss equation can be written in terms of the metric and the Fubini-Pick form as follows
\begin{align}
R_{ijkl}=&(A_{ijk,l}-A_{ijl,k})+(\chi-J)(g_{il}g_{jk}-g_{ik}g_{jl})\nnm\\ &\ +\fr2n\sum(g_{ik}A_{jlm,m}-g_{il}A_{jkm,m}) +\sum_m(A^m_{ik}A_{jlm}-A^m_{il}A_{jkm}).
\label{basic2}\end{align}

We shall use the following affine existence and uniqueness theorems later:

{\thm\label{affine existence} $($\cite{li-sim-zhao93}$)$
$($The existence$)$ Let $(M^n,g)$ be a simply connected Riemannian manifold
of dimension $n$, and $A$ be a symmetric $3$-form on $M^n$ satisfying the
affine Gauss equation \eqref{basic2} $($or equivalently \eqref{gaus}$)$ and the apolarity condition \eqref{basic1}. Then there exists a locally strongly convex immersion $x:M^n\to \bbr^{n+1}$ such that $g$ and $A$ are the affine metric and the Fubini-Pick form for $x$, respectively.}

{\thm\label{affine uniqueness} $($\cite{li-sim-zhao93}$)$ $($The uniqueness$)$ Let $x:M^n\to \bbr^{n+1}$,
$\bar x:\bar M^n\to \bbr^{n+1}$ be two locally strongly convex hypersurfaces of dimension $n$ with respectively the affine metrics $g$, $\bar g$ and the Fubini-Pick forms $A$, $\bar A$, and $\vfi:(M^n,g)\to (\bar M^n,\bar g)$ be an isometry between Riemannian manifolds. Then $\vfi^*\bar A=A$ if and only if there exists a unimodular affine transformation $\Phi:\bbr^{n+1}\to \bbr^{n+1}$ such that $\bar x\circ\vfi=\Phi\circ x$, or equivalently, $\bar x=\Phi\circ x\circ\vfi^{-1}$.}

\rmk\rm For the sufficient part of Theorem \ref{affine uniqueness}, see also \cite{lix14}.

Given a constant $L_1\in\bbr$ and a Riemannian manifold $(M^n,g)$, denote by ${\mathcal S}_{(M^n,g)}(c)$ the set of
all $TM^n$-valued symmetric bilinear forms $A\in \Gamma(\bigodot^2(T^*M^n)\bigotimes (TM^n))$, satisfying the following conditions:

(1) Under the metric $g$, the corresponding $3$-form $A\in \Gamma(\bigodot^2(T^*M^n)\bigotimes (T^*M^n))$ is totally symmetric, that is, $A\in \Gamma(\bigodot^3(T^*M^n))$;

(2) Affine Gauss equation, that is, for any $X,Y,Z\in {\mathfrak X}(M^n)$
\be\label{pre gaus_af sph1}
R(X,Y)Z=L_1(g(Y,Z)X-g(X,Z)Y)-[A(X),A(Y)](Z).
\ee

(3) $\tr_g(A)=0$,

From Theorem \ref{affine existence} and Theorem \ref{affine uniqueness}, we have
\begin{cor}\label{cor2.1}
For each $A\in {\mathcal S}_{(M^n,g)}(L_1)$, there uniquely exists one affine hypersphere $x:M^n\to\bbr^{d+1}$ with affine metric $g$, Fubini-Pick form $A$ and affine mean curvature $L_1$.
\end{cor}

Motivated by Theorem \ref{affine uniqueness}, we introduce the following concept of affine equivalence relation between nondegenerate hypersurfaces:

{\dfn Let $x:M^n\to \bbr^{n+1}$ be a nondegenerate hypersurface with the affine metric $g$. A hypersurface $\bar x:M^n\to \bbr^{n+1}$ is called affine equivalent to $x$ if there exists a unimodular transformation $\Phi:\bbr^{n+1}\to \bbr^{n+1}$ and an isometry $\vfi$ of $(M^n,g)$ such that $\bar x=\Phi\circ x\circ\vfi^{-1}$}.

To end this section, we would like to recall the following concept:

{\dfn\label{dfn afsym} {\rm(\cite{lix13})} A nondegenerate hypersurface $x:M^n\to \bbr^{n+1}$ is called affine symmetric (resp. locally affine symmetric) if

$(1)$ the pseudo-Riemannian manifold $(M^n,g)$ is symmetric (resp. locally symmetric) and therefore $(M^n,g)$ can be written (resp. locally written) as $G/K$ for some connected Lie group $G$ of isometries with $K$ one of its closed subgroups;

$(2)$ the Fubini-Pick form $A$ is invariant under the action of $G$.}

\subsection{The multiple Calabi product of hyperbolic affine hyperspheres}

For later use we make a brief review of the Calabi composition of multiple factors of hyperbolic affine hypersurfaces.
Detailed proofs of the facts listed in this subsection has been given in the articles \cite{lix11} and \cite{lix14}.
\newcommand{\stx}[2]{\strl{(#1)}{#2}}
\newcommand{\spec}[1]{\prod_{#1=1}^K\fr{c_{#1}^{n_{#1}+1}H_{(#1)}^{\fr1{n_{#1}+2}}} {(n_{#1}+1)(-\!\!\stx{#1}{L}_1)}}
\newcommand{\la}{\stx{a}{L}\!\!_1{}}\newcommand{\lb}{\stx{b}{L}\!\!_1{}}
\newcommand{\lc}{\stx{c}{L}\!\!_1{}}\newcommand{\lalp}{\stx{\alpha}{L}\!\!_1{}}
\newcommand{\ha}{\!\stx{a}{h}{}\!\!} 
\newcommand{\Ha}{H_{(a)}}\newcommand{\Hb}{H_{(b)}}\newcommand{\Hc}{H_{(c)}}
\newcommand{\ga}{\!\!\stx{a}{g}{}\!\!\!}\newcommand{\Ga}{\stx{a}{G}\!\!{}}
\newcommand{\gb}{\!\!\stx{b}{g}{}\!\!\!}\newcommand{\Gb}{\stx{b}{G}\!\!{}}
\newcommand{\galp}{\!\!\stx{\alpha}{g}{}\!\!\!}
\newcommand{\xai}{x_{a,i_a}} \newcommand{\xaij}{x_{a,i_aj_a}}
\newcommand{\HH}{f_K\prod_a\fr{c_a^{(n_a+1)(f_K-1)}\Ha^{\fr{f_K+1}{n_a+2}}}
{(n_a+1)^{f_K-n_a}(-\!\!\la)^{f_K-n_a-1}}}
\newcommand{\h}{f^{-\fr1{n+2}}_K \prod_a\fr {(n_a+1)^{\fr{f_K-n_a}{f_K+1}}(-\!\!\la)^{\fr{f_K-n_a-1}{f_K+1}}}
{\left(c_a^{n_a+1}\right)^{\fr{f_K-1}{f_K+1}}\Ha^{\fr1{n_a+2}}}}
\newcommand{\oma}{\stx{a}{\omega}{}\!\!}
\newcommand{\omb}{\stx{b}{\omega}{}\!\!}
\newcommand{\tdca}{(n_a+1)\big(-\!\!\la\big)\prod_b\fr{c_b^{n_b+1}} {(n_b+1)\big(-\!\!\lb\big)}}
\newcommand{\cha}{\prod_a\fr{c_a^{n_a+1}\Ha^{\fr1{n_a+2}}} {(n_a+1)\big(-\!\!\la\big)}}
\newcommand{\chb}{\prod_b\fr{c_b^{n_b+1}\Hb^{\fr1{n_b+2}}} {(n_b+1)\big(-\!\!\lb\big)}}
\newcommand{\chc}{\prod_c\fr{c_c^{n_c+1}\Hc^{\fr1{n_c+2}}} {(n_c+1)\big(-\!\!\lc\big)}}
\newcommand{\ca}{\prod_a\fr{c_a^{n_a+1}}{(n_a+1)\big(-\!\!\la)}}
\newcommand{\cb}{\prod_b\fr{c_b^{n_b+1}}{(n_b+1)\big(-\!\!\lb)}}
\newcommand{\cc}{\prod_c\fr{c_c^{n_c+1}}{(n_c+1)\big(-\!\!\lc)}}
\newcommand{\Aa}{\stx{a}{A}{}\!\!}
\newcommand{\Aalp}{\stx{\alpha}{A}{}\!\!}
\newcommand{\olomea}{\stx{a}{\ol\omega}{}\!\!}
\newcommand{\olomeb}{\stx{b}{\ol\omega}{}\!\!}
\newcommand{\olgma}{\stx{a}{\ol\Gamma}{}\!\!\!}
\newcommand{\olgmb}{\stx{b}{\ol\Gamma}{}\!\!\!}

Let $r,s$ be two nonnegative integers with $K:=r+s\geq 2$ and $x_\alpha:M^{n_\alpha}_\alpha\to\bbr^{n_\alpha+1}$, $1\leq \alpha\leq s$, be hyperbolic affine hyperspheres of dimension $n_\alpha>0$ with affine mean curvatures $\stx{\alpha}{L}\!\!_1$ and with the origin their common affine center. For convenience we make the following convention:
$$1\leq a,b,c\cdots\leq K,\quad 1\leq\lambda,\mu,\nu\leq K-1,\quad
1\leq\alpha,\beta,\gamma\leq s,\quad \td\alpha=\alpha+r,\ \td\beta=\beta+r,\ \td\gamma=\gamma+r.
$$
Furthermore, for each $\alpha=1,\cdots,s$, set $\td i_{\alpha}=i_\alpha+K-1+\sum_{\beta<\alpha}n_\beta$ with $1\leq i_\alpha\leq n_\alpha$.

Define
$$
f_a:=\begin{cases} a,&1\leq a\leq r;\\ \sum_{\beta\leq \alpha}n_\beta+\td{\alpha},&r+1\leq a=\td\alpha\leq r+s,
\end{cases}
$$
and
$$e_a:=\exp\left(-\fr{t_{a-1}}{n_{a}+1}+\fr{t_{a}}{f_{a}}+\fr{t_{a+1}}{f_{a+1}} +\cdots+\fr{t_{K-1}}{f_{K-1}}\right),\quad 1\leq a\leq K=r+s$$
In particular,
$$
e_1=\exp\left(\fr{t_1}{f_1}+\fr{t_2}{f_2} +\cdots+\fr{t_{K-1}}{f_{K-1}}\right),\quad
e_K=\exp\left(-\fr{t_{K-1}}{n_K+1}\right).
$$

Put $n=\sum_\alpha n_\alpha+K-1$ and $M^n=R^{K-1}\times M^{n_1}_1\times\cdots\times M^{n_s}_s$. For any $K$ positive numbers $c_1,\cdots,c_K$, define a smooth map $x:M^n\to\bbr^{n+1}$ by
\begin{align}
x(t^1,&\cdots,t^{K-1},p_1,\cdots,p_s):=(c_1e_1,\cdots, c_re_r,c_{r+1}e_{r+1}x_1(p_1),\cdots,c_Ke_Kx_s(p_s)),\nnm\\&\hs{1cm}\forall (t^1,\cdots,t^{K-1},p_1,\cdots,p_s)\in M^n.\label{mulpro2}
\end{align}

{\prop\label{general sense} (\cite{lix11}) The map $x:M^n\to\bbr^{n+1}$ defined above is a new hyperbolic affine hypersphere with the affine mean curvature
\be\label{newl1c}
L_1=-\fr1{(n+1)C},\quad C:=\left(\fr1{n+1}\prod_{a=1}^r c_a^2\cdot\prod_{\alpha=1}^s\fr{c_{r+\alpha}^{2(n_\alpha+1)}} {(n_\alpha+1)^{n_\alpha+1}(-\!\!\stx{\alpha}{L}_1)^{n_\alpha+2}}\right)^{\fr1{n+2}},
\ee
Moreover, for given positive numbers $c_1,\cdots,c_K$, there exits some $c>0$ and $c'>0$ such that
the following three hyperbolic affine hyperspheres
\bea &x:=(c_1e_1,\cdots, c_re_r,c_{r+1}e_{r+1}x_1,\cdots,c_Ke_sx_s),\nnm\\
&\bar x:=c(e_1,\cdots, e_r,e_{r+1}x_1,\cdots,e_sx_s),\nnm\\
&\td x:=(e_1,\cdots, e_r,e_{r+1}x_1,\cdots,c'e_sx_s)\nnm
\eea
are equiaffine equivalent to each other.}

{\dfn\label{df2} {\rm(\cite{lix11})} \rm The hyperbolic affine hypersphere $x$ is called the Calabi composition of $r$ points and $s$ hyperbolic affine hyperspheres.}

Then we have

{\cor\label{cor} {\rm(\cite{lix11})} The Calabi composition $x:M^n\to \bbr^{n+1}$ of $r$ points and $s$ hyperbolic affine hyperspheres $x_\alpha:M^{n_\alpha}\to\bbr^{n_\alpha+1}$, $1\leq \alpha\leq s$, is affine symmetric if and only if
each positive dimensional factor $x_\alpha$ is symmetric.}

Note that for a given locally strongly convex hypersurface $x:M^n\to\bbr^{n+1}$ with the affine metric $g$, $(M^n,g)$ is a Riemannian manifold. Then we have the following characterization of Calabi composition of symmetric factors which is important in the proof of Theorem \ref{main}:

{\thm\label{chara} {\rm(\cite{lix14}; cf. \cite{lix13})} A locally strongly convex and affine symmetric {\bf hypersphere} $x:M^n\to\bbr^{n+1}$ is locally affine equivalent to the Calabi composition of some hyperbolic affine hyperspheres possibly including point factors if and only if $M^n$ is reducible as a Riemannian manifold with respect to the affine metric.}

\section{Some typical examples}

To make the main theorem more understandable, we provide in this section a systematic and unified treatment of some typical examples of affine symmetric hyperspheres in $\bbr^{n+1}$ giving, for the first time, the necessary computation details. These examples have partly appeared in \cite{sas80}, \cite{li-sim-zhao93}, \cite{nom-sas94}, \cite{dil-vra94}, \cite{bir-djo12}, \cite{lix13} and particularly in the important classification theorem by Z.J. Hu, H.Z. Li and L. Vrancken (\cite{hu-li-vra11}, see also Theorem \ref{cla thm} in the next section).

\expl\label{expl1} \rm (\cite{li-sim-zhao93}, \cite{nom-sas94}) Quadric Hypersurfaces

There are three kinds of quadric hypersurfaces in $\bbr^{n+1}$ and they are given by the following quadric equations
\begin{align}
&\text{(1)\ Ellipsoid:}\quad(x^1)^2+\cdots+(x^n)^2+(x^{n+1})^2=c^2,\quad c>0;\\
&\text{(2)\ Paraboloid:}\quad(x^1)^2+\cdots+(x^n)^2=2x^{n+1};\\
&\text{(3)\ Hyperboloid:}\quad(x^1)^2+\cdots+(x^n)^2-(x^{n+1})^2=-c^2,\quad x^{n+1}>0,
\quad c>0.\hs{1.1in} \end{align}

It is well known that the above three hypersurfaces are (resp. elliptic, hyperbolic and parabolic) affine hyperspheres (with resp. positive, negative and zero affine principal curvatures) and have vanishing Fubini-Pick forms. It then follows that, with respect to the affine metrics, they have constant (resp. positive, negative and zero) affine sectional curvatures. In particular, they are affine symmetric hyperspheres. Also we have

{\prop\label{expl1-prop} {\rm(\cite{li-sim-zhao93})} A locally strongly convex hypersurface $x:M\to\bbr^{n+1}$ has vanishing Fubini-Pick form if and only if it is one of the above quadric hypersurfaces.}

\expl\label{expl2} \rm (\cite{li-sim-zhao93}) The standard flat hypersurfaces with nonzero Fubini-Pick form

Given a positive number $C$, let $x:\bbr^{n}\to \bbr^{n+1}$ be the well known flat hyperbolic affine hypersphere of dimension $n$ which is defined by
$$
x^1\cdots x^{n} x^{n+1}=C,\quad x^1>0,\cdots,x^{n+1}>0.
$$
Then it is not hard to see that $x$ is the Calabi composition of $n+1$ points and thus is affine flat. In fact, we can write for example
$$
x=(e_1,\cdots,e_{n},Ce_{n+1}).
$$
It follows from Corollary \ref{cor} that $x$ is affine symmetric. In particular, $x$ has a positive constant Pick invariant.

Note that by a theorem of L. Vrancken, A-M. Li and U. Simon in \cite{vra-li-sim91} (also see \cite{amli89}), Example \ref{expl2} is, up to equiaffine equivalence, the only one with flat affine metric and positive Pick invariant.

\expl\label{expl3} \rm(\cite{nom-sas94}, \cite{dil-vra94}, \cite{hu-li-sim-vra09}) The standard embedding
$$x:M\equiv{\rm SL}(m,\bbr)/{\rm SO}(m)\to\bbr^{n+1},\quad n=\fr12m(m+1)-1,\quad m\geq 3.$$

Let $\mathfrak{s}\mathfrak{l}(m,\bbr)$, $\mathfrak{s}\mathfrak{o}(m)$ be the Lie algebras of ${\rm SL}(m,\bbr)$, ${\rm SO}(m)$ respectively, and $\bbr^{n+1}\equiv {\mathfrak s}(m)$ the vector space of real symmetric matrices of order $m$. Then the canonical decomposition of $\mathfrak{s}\mathfrak{l}(m,\bbr)$ with respective to $\mathfrak{s}\mathfrak{o}(m)$ is ${\mathfrak s}{\mathfrak l}(m,\bbr)={\mathfrak s}{\mathfrak o}(m,\bbr)+{\mathfrak s}_0(m)$ where
$${\mathfrak s}_0(m):=\{X\in {\mathfrak s}(m);\ \tr X=0\}$$
and is naturally identified with the tangent space $T_oM$ at the origin $o={\rm SO}(m)\in M$, the coset of the identity matrix.

There is a representation $\phi$ of ${\rm SL}(m,\bbr)$ on $\bbr^{n+1}$ defined by
$$
\phi(a)X:=aXa^t,\quad \text{for\ }a\in {\rm SL}(m,\bbr),\ X\in \bbr^{n+1}.
$$
Then we have

{\lem\label{expl3-lem}{\rm(\cite{nom-sas94})} $\phi({\rm SL}(m,\bbr))\subset {\rm SL}(n+1,\bbr)$. So $\phi({\rm SL}(m,\bbr))$ can be taken to be a subgroup of the unimodular group ${\rm UA}(n+1)$ on $\bbr^{n+1}$.}

For a given constant $L_1<0$, put
$$C=\fr{\sqrt{m}}{4}\left(\fr4{m(-L_1)}\right)^{\fr{n+2}2}$$
and define a map
$x:{\rm SL}(m,\bbr)/{\rm SO}(m)\to \bbr^{n+1}$ as follows:
$$
x(\,{}_a{\rm SO}(m))=Caa^t,\quad \text{for\ }a\in {\rm SL}(m,\bbr).
$$
Then it is clear that $x$ is equivariant with respect to the representation $\phi:{\rm SL}(m,\bbr)\to {\rm UA}(n+1)$ (see Lemma \ref{expl3-lem}) and $x(M)$ coincides with the subset of all positive-definite matrices in ${\mathfrak s}(m)$ with constant determinant $C^m$, and $x(o)=CI_m$ where $I_m$ is the identity matrix of order $m$.

Furthermore, $x$ is an equiaffine symmetric hypersphere of affine mean curvature $L_1$. In fact, this last conclusion follows by the following computation:

Now for each $X\in{\mathfrak s}_0(m)\equiv T_oM$, $a(t):={}_{\exp tX}{\rm SO}(m)$ is a geodesic curve on $M$. Then it holds that
$$
x_*(X)=\left.\dd{}{t}\right|_{t=0}x(a(t))=C\left.\dd{}{t}\right|_{t=0}((\exp tX)(\exp tX)^t)=2CX.
$$
This shows that $x$ is an immersion at $o$ and thus is an immersion globally since $x$ is equivariant. Clearly, $x$ is injective and is thus an imbedding of $M$ into $\bbr^{n+1}$.

Moreover, the standard inner product $(\cdot,\cdot)$ on $\bbr^{n+1}\equiv {\mathfrak s}(m)$ is defined by $(X,Y)=\tr(XY)$, $X,Y\in {\mathfrak s}(m)$.
Since
$$(x_*(X),x(o))=(2CX,CI_m)=2C^2\tr(XI_m)=2C^2\tr X=0,\quad X\in T_oM\equiv{\mathfrak s}_0(m),$$
$x(o)$ is a transversal vector of $x$ at $o$ and thus is transversal everywhere by the equivariance.

On the other hand, if we denote by $Y^*$ the Killing vector field on $M$ induced by $Y\in {\mathfrak s}_0(m)$, then the value of $Y^*$ at $a(t)$
$$
Y^*|_{a(t)}=\left.\dd{}{s}\right|_{s=0}(\,{}_{\exp sY a(t)}{\rm SO}(m))=\left.\dd{}{s}\right|_{s=0}({}\,_{\exp s Y \exp t X}{\rm SO}(m)).
$$
Therefore
$$
x_*(Y^*|_{a(t)})=C\left.\dd{}{s}\right|_{s=0}((\exp sY \exp tX)(\exp sY \exp tX)^t).
$$
It follows that
\begin{align}
X(x_*(Y^*))=&C\left.\ppp{}{t}{s}\right|_{t=s=0}(\exp sY \exp tX \exp sX^t \exp tY^t)\nnm\\
=&2C(YX+XY)=2C\left(YX+XY-\fr2m\tr(XY)I_m\right)+\fr4mC(X,Y)I_m\label{expl3-gausf}
\end{align}
implying that $x$ is locally strongly convex since $(X,Y)$ is positive definite.
Moreover the affine metric (Blaschke metric) of $x$ at the origin $o$ is by definition
$$
g_o(X,Y)=\left(\fr{4C}{\sqrt{m}}\right)^{\fr2{n+2}}(X,Y)=-\fr4{mL_1}(X,Y),\quad X,Y\in {\mathfrak s}_0(m).
$$
Clearly $g_o$ is positive definite and invariant by ${\rm SO}(m)$ and it induced a invariant Riemannian metric $g$. On the other hand, the involution map $\sigma:{\mathfrak s}{\mathfrak l}(m,\bbr)\to {\mathfrak s}{\mathfrak l}(m,\bbr)$ defined by $\sigma(X)=-X^t$ is isometric with respect to $g_o$, thus the invariant metric $g$ is symmetric; Note that $x$ is equivariant, thus $g$ is nothing but the affine metric of $x$.

Let $A_o$ be the $(1,2)$ tensor on $\mathfrak{s}_0(m)$ defined by
$$
A_o(X,Y)=XY+YX-\fr2m\tr(XY)I_m,\quad \forall X,Y\in\mathfrak{s}_0(m),
$$
which gives a linear map for any $X\in \mathfrak{s}_0(m)$: $A_o(X):\mathfrak{s}_0(m)\to \mathfrak{s}_0(m)$ by $A_o(X)Y=A_o(X,Y)$, $Y\in \mathfrak{s}_0(m)$.

To find the affine normal vector at $o$, we should first prove the following lemma:

{\lem\label{expl3-lem2} Define $A_o(X,Y,Z)=g_o(A_o(X,Y),Z)$, for $X,Y,Z\in \mathfrak{s}_0(m)$. Then

$(1)$ the $(0,3)$-tensor $A_o(X,Y,Z)$ is totally symmetric;

$(2)$ for each $X\in \mathfrak{s}_0(m)$, the linear map $A_o(X)$ is traceless.}

\proof
Conclusion (1) is direct. To prove (2), we denote by $e^j_i$ the $m\times m$ matrix with the $(i,j)$-th element being $1$ and all other elements zero, that is, its $(k,l)$-th element $(e^j_i)^k_l=\delta^k_i\delta^j_l$, $1\leq k,l\leq m$. Then $\{e^j_i,\ 1\leq i,j\leq m\}$ is the standard basis for the real linear space $M(m,\bbr)$ of $m\times m$ real matrices. Define
$$
f_\alpha=e^\alpha_\alpha-e^m_m\text{\ for\ }1\leq \alpha\leq m-1;\quad f^j_i=\fr12(e^j_i+e^i_j)\text{\ for\ }1\leq i<j\leq m.
$$
Then $\{f_\alpha,f^j_i\}$ is a basis for $\mathfrak{s}_0(m)$. For $X=(X^i_j)\in \mathfrak{s}_0(m)$, we find by direct computation
\begin{align}
A_o(X)f_\alpha=&f_\alpha X+Xf_\alpha-\fr2m\tr(f_\alpha X)I_m=\fr2m((m-1)X^\alpha_\alpha+X^m_m)f_\alpha+\cdots,\label{expl3-21}\\
A_o(X)f^j_i=&f^j_iX+Xf^j_i-\fr2m\tr(f^j_iX)I_m=(X^j_j+X^i_i)f^j_i+\cdots\label{expl3-22}
\end{align}
where we have omitted those terms not containing $f_\alpha$ in \eqref{expl3-21}, and those not containing $f^j_i$ in \eqref{expl3-22}, respectively. It then follows that
\begin{align}
\tr A_o(X)=&\fr2m\sum_\alpha((m-1)X^\alpha_\alpha+X^m_m) +\sum_{i<j}(X^j_j+X^i_i)\nnm\\ =&\fr{2(m-1)}m\sum_iX^i_i+\fr12\sum_{i,j}(X^i_i+X^j_j)-\sum_iX^i_i=0
\end{align}
since $\tr X=\sum_iX^i_i=0$.
\endproof

Since $Y^*$ is chosen to be the Killing vector field on $M$ corresponding to $Y$, we have $\hat\nabla_XY^*=0$ where $\hat\nabla$ is the Levi-Civita connection of the affine metric $g$. Therefore by taking the trace of \eqref{expl3-gausf} with respect to $g_o$ and using Lemma \ref{expl3-lem2}, we find that, at $o$, the affine normal vector
$$
\xi_o=\fr1n\Delta_g x=\left(\fr{4C}{\sqrt{m}}\right)^{-\fr2{n+2}}\fr4m\cdot x(o)=-L_1x(o).
$$
$\xi_o$ is clearly invariant by $\phi({\rm SO}(m))$ and for any $X\in {\mathfrak s}{\mathfrak l}(m,\bbr)$, $\phi_*(X)\xi_o\in x_*(T_o(M))$. Then the equivariant transversal vector field $\xi$ induced by $\xi_o$ coincides with the affine normal vector (see Lemma 4.4 in \cite{nom-sas94}). Since $x$ is also equivariant, $\xi=-L_1x$ holds identically. Therefore, $x$ is a hyperbolic affine sphere with affine mean curvature $L_1$.

Now the equivariance of $x$ implies that its Fubini-Pick form $A$ is ${\rm SL}(m,\bbr)$-invariant which indicates that $x$ is affine symmetric, and the invariant Fubini-Pick form $A$ is uniquely determined by the cubic form $A_o$ given in Lemma \ref{expl3-lem2} (see also Definition \eqref{f-p}):
\be\label{expl3-eq5}
A_o(X,Y,Z)=g_0\left(\left(XY+YX-\fr2m\tr(XY)I_m\right),Z\right),\quad X,Y,Z\in \mathfrak{s}_0(m).
\ee

\expl\label{expl4} \rm(\cite{hu-li-vra11}, cf. \cite{bir-djo12} for $m=3$) The standard embedding
$$x:M\equiv{\rm SL}(m,\bbc)/{\rm SU}(m)\to\bbr^{n+1}\quad n=m^2-1,\quad m\geq 3.$$

Let $\mathfrak{s}\mathfrak{l}(m,\bbc)$, $\mathfrak{s}\mathfrak{u}(m)$ be the Lie algebras of ${\rm SL}(m,\bbc)$, ${\rm SU}(m)$ respectively, and $\bbr^{n+1}\equiv {\mathfrak h}(m)$ the vector space of complex Hermitian matrices of order $m$. Then the canonical decomposition of $\mathfrak{s}\mathfrak{l}(m,\bbc)$ with respective to $\mathfrak{s}\mathfrak{u}(m)$ is ${\mathfrak s}{\mathfrak l}(m,\bbc)={\mathfrak s}{\mathfrak u}(m)+{\mathfrak h}_0(m)$ where
$${\mathfrak h}_0(m):=\{X\in {\mathfrak h}(m);\ \tr X=0\},$$
which can be identified with the tangent space $T_oM$ at the origin $o={\rm SU}(m)\in M$.

There is a representation $\phi$ of ${\rm SL}(m,\bbc)$ on $\bbr^{n+1}$ by
$$
\phi(a)X:=aX\bar a^t,\quad \text{for\ }a\in {\rm SL}(m,\bbc),\ X\in \bbr^{n+1}.
$$

{\lem\label{expl4-lem}\rm $\phi({\rm SL}(m,\bbc))\subset {\rm SL}(n+1,\bbr)$ and thus $\phi({\rm SL}(m,\bbc))$ can be viewed as a subgroup of the unimodular group ${\rm UA}(n+1)$ on $\bbr^{n+1}$.}

\proof Let $e^j_i$ and $f^j_i$ be as in Example \ref{expl3}. Then $\{e^j_i,\ 1\leq i,j\leq m\}$ can also be taken as the standard basis for the complex linear space $M(m,\bbc)$ of $m\times m$ complex matrices, with its complex dual basis denoted by $\{\omega^i_j,\ 1\leq i,j\leq m\}$. Define
$$
\td f^j_i=\fr12\sqrt{-1}(e^j_i-e^i_j)\text{\ for\ }1\leq i<j\leq m.
$$
Then $\{e^i_i,f^j_i,\td f^j_i\}$ is a basis for the real linear space $\mathfrak{h}(m)$ with the dual basis $\{\theta^i_i,\theta^i_j,\td \theta^i_j\}$ where
$$
\theta^i_i=\omega^i_i\text{\ for\ }1\leq i\leq m;\quad \theta^i_j=(\omega^i_j+\omega^j_i),\ \td \theta^i_j=\sqrt{-1}(\omega^j_i-\omega^i_j)\text{\ for\ }1\leq i<j\leq m.
$$
It follows that for $i<j$,
\begin{align}
\theta^i_j(e^l_k)=&\omega^i_j(e^l_k)+\omega^j_i(e^l_k) =\delta^i_k\delta^l_j+\delta^j_k\delta^l_i,\label{expl4-6}\\
\td\theta^i_j(e^l_k)=&\sqrt{-1}(\omega^j_i(e^l_k)-\omega^i_j(e^l_k)) =\sqrt{-1}(\delta^j_k\delta^l_i-\delta^i_k\delta^l_j).\label{expl4-7}
\end{align}

For each $X\in\mathfrak{s}\mathfrak{l}(m,\bbc)$, write $X=(X^k_l)_{m\times m}=\sum_{k,l}X^k_le^l_k$. Then $\tr X=\sum_iX^i_i=0$ and, for each pair of $i,j$, we have
\be\label{expl401}
Xe^j_i=\sum_{k,l,p}(X^k_p\delta^p_i\delta^j_l)e^l_k=\sum_{k}X^k_ie^j_k,\quad e^j_i\bar X^t=\sum_{k,l,p}(\delta^k_i\delta^j_p\bar X^l_p)e^l_k=\sum_{k}\bar X^k_je^k_i.
\ee
Since, by definition, $\phi_*(X)(A)=XA+A\bar X^t$ ($X\in {\mathfrak s}{\mathfrak l}(m,\bbc)$, $A\in\bbr^{n+1}$), it follows by \eqref{expl4-6}--\eqref{expl401} and $\tr X=0$ that
\begin{align}
\sum_i\theta^i_i(\phi_*(X)(f^i_i))=&\sum_i\omega^i_i(Xe^i_i+e^i_i\bar X^t) =\sum_{i,k}(X^k_i\omega^i_i(e^i_k)+\bar X^k_i\omega^i_i(e^k_i))\nnm\\
=&\sum_i(X^i_i+\bar X^i_i)=0,\label{expl4-8}\\
\sum_{i<j}\theta^i_j(\phi_*(X)(f^j_i)) =&\fr12\sum_{i<j}\theta^i_j(X(e^j_i+e^i_j)+(e^j_i+e^i_j)\bar X^t)\nnm\\
=&\fr12\sum_{i<j,k}(X^k_i\theta^i_j(e^j_k)+X^k_j\theta^i_j(e^i_k)+\bar X^k_j\theta^i_j(e^k_i)+\bar X^k_i\theta^i_j(e^k_j))\nnm\\
=&\fr12\sum_{i<j}(X^i_i+X^j_j+\bar X^j_j+\bar X^i_i)=\fr12\sum_{i\neq j}(X^i_i+\bar X^i_i)\nnm\\
=&\fr{m-1}2\sum_i(X^i_i+\bar X^i_i)=0,\label{expl4-9}\\
\sum_{i<j}\td\theta^i_j(\phi_*(X)(\td f^j_i))=&\fr12\sqrt{-1}\sum_{i<j}\td\theta^i_j(X(e^j_i-e^i_j)+(e^j_i-e^i_j)\bar X^t) \nnm\\ =&\fr12\sqrt{-1}\sum_{i<j,k}(X^k_i\td\theta^i_j(e^j_k)-X^k_j\td\theta^i_j(e^i_k) +\bar X^k_j\td\theta^i_j(e^k_i)-\bar X^k_i\td\theta^i_j(e^k_j))\nnm\\
=&\fr12\sum_{i<j}(X^i_i+X^j_j+\bar X^j_j+\bar X^i_i) =\fr12\sum_{i\neq j}(X^i_i+\bar X^i_i)\nnm\\
=&\fr{m-1}2\sum_i(X^i_i+\bar X^i_i)=0.\label{expl4-10}
\end{align}

Taking the sum of \eqref{expl4-8}--\eqref{expl4-10}, we find
$$
\tr(\phi_*(X))=\sum_i\theta^i_i(\phi_*(X)(f^i_i)) +\sum_{i<j}\theta^i_j(\phi_*(X)(f^j_i)) +\sum_{i<j}\td\theta^i_j(\phi_*(X)(\td f^j_i))=0,
$$
completing the proof of Lemma \ref{expl4-lem}.\endproof

For a given constant $L_1<0$, put
$$C=\fr{\sqrt{m}}{4}\left(\fr4{m(-L_1)}\right)^{\fr{n+2}2}$$
and define a map
$x:{\rm SL}(m,\bbc)/{\rm SU}(m)\to \bbr^{n+1}$ as follows:
$$
x(\,{}_a{\rm SU}(m))=Ca\bar a^t,\quad \text{for\ }a\in {\rm SL}(m,\bbc).
$$
Then, by Lemma \ref{expl4-lem}, the $x$ is equivariant with respect to the representation $\phi:{\rm SL}(m,\bbc)\to {\rm UA}(n+1)$.

Now for each $X\in \mathfrak{h}_0(m)$, define $a(t)={}_{\exp tX}{\rm SO}(m)$. Then
$$
x_*(X)=\left.\dd{}{t}\right|_{t=0}x(a(t))=C\left.\dd{}{t}\right|_{t=0}((\exp tX)(\ol{\exp tX})^t)=2CX.
$$
Thus $x$ is an immersion at $o$ and thus everywhere. Moreover, $x$ is also an imbedding of $M$ into $\bbr^{n+1}$.

For $X,Y\in {\mathfrak h}(m)$, define $(X,Y)=\tr(XY)$. Then $(\cdot,\cdot)$ is the standard inner product on $\bbr^{n+1}\equiv {\mathfrak h}(m)$. In particular, it is positive definite. As in Example \ref{expl3}, $x$ is equivariant and transversal everywhere on $M$.

For any $Y\in \mathfrak{h}_0(m)$, the corresponding Killing vector field $Y^*$ on $M$ satisfies
$$
Y^*|_{a(t)}=\left.\dd{}{s}\right|_{s=0}(\,{}_{\exp sY a(t)}{\rm SU}(m))=\left.\dd{}{s}\right|_{s=0}({}\,_{\exp sY \exp tX}{\rm SU}(m)).
$$
It then follows that
$$
x_*(Y^*|_{a(t)})=C\left.\dd{}{s}\right|_{s=0}((\exp sY \exp tX)(\ol{\exp sY \exp tX})^t).
$$
Therefore
\begin{align}
X(x_*(Y^*))=&C\left.\ppp{}{t}{s}\right|_{t=s=0}(\exp sY \exp tX \exp s\ol X^t \exp t\ol Y^t)\nnm\\
=&2C(XY+YX)=2C\left(YX+XY-\fr2m\tr(XY)I_m\right)+\fr4mC(X,Y)I_m\label{expl4-gausf}
\end{align}
implying that $x$ is locally strongly convex as $(X,Y)$ is positive definite.
Thus, at the origin $o$, the invariant affine metric $g_o$ is defined by:
$$
g_o(X,Y)=\left(\fr{4C}{\sqrt{m}}\right)^{\fr2{n+2}}(X,Y)=-\fr4{mL_1}(X,Y),\quad X,Y\in \mathfrak{h}_0(m).
$$
Since $g_o$ is positive definite and invariant by ${\rm SU}(m)$, the invariant Riemannian metric $g$ determined by $g_o$ is exactly the affine metric of $x$. Similar to Example \ref{expl3}, we can prove that, for any $X\in \mathfrak{h}_0(m)$, the real linear map
$$
Y\in \mathfrak{h}_0(m)\mapsto XY+YX-\fr2m\tr(XY)I_m
$$
is also traceless. So, by making use of \eqref{expl4-gausf}, we find that $\xi=-L_1x$ holds identically, implying that $x$ is a hyperbolic affine sphere with affine mean curvature $L_1$.

Note that the involution map $\sigma:{\mathfrak s}{\mathfrak l}(m,\bbc)\to {\mathfrak s}{\mathfrak l}(m,\bbc)$ is given by $\sigma(X)=-\bar X^t$ and isometric with respect to $g_o$, thus the invariant affine metric $g$ is symmetric; Furthermore, the Killing vector field $Y^*$ on $M$ given by $Y\in \mathfrak{h}_0(m)$ subject to $\hat\nabla_XY^*=0$ for all $X\in \mathfrak{h}_0(m)$ with $\hat\nabla$ the Levi-Civita connection. It follows from \eqref{expl3-gausf} and Definition \eqref{f-p} that the Fubini-Pick form $A$ of $x$ is invariant and is determined by its value $A_o$ at the origin $o$:
\be\label{expl4-eq5}
A_o(X,Y,Z)=g_0\left(\left(XY+YX-\fr2m\tr(XY)I_m\right),Z\right),\quad X,Y,Z\in \mathfrak{h}_0(m).
\ee
This shows that $x$ is an affine symmetric hypersphere.

\expl\label{expl5}\rm(\cite{hu-li-vra11}, cf. \cite{bir-djo12} for $m=3$) The standard embedding
$$x:M\equiv{\rm SU}^*(2m)/{\rm Sp}(m)\to\bbr^{n+1}\quad n=2m^2-m-1,\quad m\geq 3,$$
where ${\rm SU}^*(2m)={\rm SL}(2m,\bbc)\cap {\rm U}^*(2m)$ with ${\rm U}^*(2m)$ the usual ${\rm U}$-star group of order $2m$.

Define $J=\lmx 0&-I_m\\I_m&0\rmx$. Then the $U$-star group, or in other words, the general quaternion linear group has an expression in terms of complex matrices as
\begin{align}
U^*(2m)=&\{T\in GL(2m,\bbc);\ TJ=J\bar T\}\nnm\\
=&\{T=\lmx A&B\\-\bar B&\bar A\rmx\in {\rm GL}(2m,\bbc);\ A,B\in {\rm M}(m,\bbc)\}.
\end{align}

Consequently, the Lie algebra of $U^*(2m)$ is written as
\begin{align}
\mathfrak{u}^*(2m)=&\left\{X\in {\rm M}(2m,\bbc);\ XJ=J\bar X\right\}\nnm\\
=&\{X=\lmx A&B\\-\bar B&\bar A\rmx\in {\rm M}(2m,\bbc);\ A,B\in {\rm M}(m,\bbc)\}.
\end{align}
It follows that the special $U$-star group or the special quaternion linear group ${\rm SU}^*(2m)$ is given by
$$
{\rm SU}^*(2m)={\rm SL}(2m,\bbc)\cap {\rm U}^*(2m)=\left\{T\in U^*(2m);\ \det T=1\right\}
$$
of which the Lie algebra is
$$
\mathfrak{s}\mathfrak{u}^*(2m)=\left\{X\in \mathfrak{u}^*(2m),\ \tr X=0.\right\}
$$
Moreover, the quaternion unitary group or the symplectic group is defined by
$${\rm Sp}(m)={\rm U}(2m)\cap {\rm SU}^*(2m)=\{T\in {\rm SU}^*(2m);\ T\bar T^t=I_{2m}\}$$ with the Lie algebra
$$\mathfrak{s}\mathfrak{p}(m)=\left\{X\in\mathfrak{s}\mathfrak{u}^*(2m);\ X+\bar X^t=0.\right\}.
$$

Let $\bbr^{n+1}\equiv {\mathfrak q}{\mathfrak h}(m)$ be the real vector space of quaternion Hermitian matrices of order $m$. Then we have
$$
{\mathfrak q}{\mathfrak h}(m)=\mathfrak{h}(m)\oplus \mathfrak{s}\mathfrak{o}(m,\bbc)
=\{\lmx A&B\\-\bar B&\bar A\rmx\in {\rm M}(2m,\bbc);\ A\in \mathfrak{h}(m),B\in\mathfrak{s}\mathfrak{o}(m,\bbc)\}
$$
there is a representation $\phi$ of ${\rm SU}^*(2m)$ on $\bbr^{n+1}$ by
$$
\phi(a)X:=aX\bar a^t,\quad \text{for\ }a\in {\rm SU}^*(2m),\ X\in \bbr^{n+1}.
$$
Suitably choose a basis for the real vector space $\bbr^{n+1}$ together with its dual basis, and then by a similar computation as in Example \ref{expl4} we are able to obtain

{\lem\label{expl5-lem}\rm $\phi({\rm SU}^*(2m))\subset {\rm SL}(n+1,\bbr)$, that is, $\phi({\rm SU}^*(2m))$ can be viewed as a subgroup of the unimodular group ${\rm UA}(n+1)$ on $\bbr^{n+1}$.}

Define ${\mathfrak q}{\mathfrak h}_0(m)=\{X\in{\mathfrak q}{\mathfrak h}(m);\ \tr X=0\}$.
Then the canonical decomposition of $\mathfrak{s}\mathfrak{u}^*(2m)$ with respect to $\mathfrak{s}\mathfrak{p}(m)$ is as follows:
$$
\mathfrak{s}\mathfrak{u}^*(2m)=\mathfrak{s}\mathfrak{p}(m)
+{\mathfrak q}{\mathfrak h}_0(m)
$$
where the subspace ${\mathfrak q}{\mathfrak h}_0(m)$ can be identified with the tangent space $T_oM$ at the origin $o={\rm Sp}(m)\in M$.

For a given constant $L_1<0$, put
$$C=\fr{\sqrt{2m}}{4}\left(\fr2{m(-L_1)}\right)^{\fr{n+2}2}$$
and define a map
$x:{\rm SU}^*(2m)/{\rm Sp}(m)\to \bbr^{n+1}$ as follows:
$$
x(\,{}_a{\rm Sp}(m))=Ca\bar a^t,\quad \text{for\ }a\in {\rm SU}^*(m).
$$
Then, by Lemma \ref{expl5-lem}, the $x$ is equivariant with respect to the representation $\phi:{\rm SU}^*(2m)\to {\rm UA}(n+1)$.

Now for each $X\in {\mathfrak q}{\mathfrak h}_0(m)$, define $a(t)={}_{\exp t X}{\rm Sp}(m)$. Then
$$
x_*(X)=\left.\dd{}{t}\right|_{t=0}x(a(t))=C\left.\dd{}{t}\right|_{t=0}((\exp tX)(\ol{\exp tX})^t)=2CX.
$$
Thus $x$ is an immersion at $o$ and thus everywhere. Moreover, $x$ is also an imbedding of $M$ into $\bbr^{n+1}$.

For $X,Y\in {\mathfrak q}{\mathfrak h}(m)$, define $(X,Y)=\tr(XY)$. Then $(\cdot,\cdot)$ is the standard inner product on $\bbr^{n+1}\equiv {\mathfrak q}{\mathfrak h}(m)$. In particular, it is positive definite. As in Example \ref{expl3}, $x$ is invariant and transversal everywhere on $M$.

For any $Y\in {\mathfrak q}{\mathfrak h}_0(m)$, the corresponding Killing vector field $Y^*$ on $M$ satisfies
$$
Y^*|_{a(t)}=\left.\dd{}{s}\right|_{s=0}(\,{}_{\exp sY a(t)}{\rm Sp}(m))=\left.\dd{}{s}\right|_{s=0}({}\,_{\exp sY \exp tX}{\rm Sp}(m)).
$$
It then follows that
$$
x_*(Y^*|_{a(t)})=C\left.\dd{}{s}\right|_{s=0}((\exp sY \exp tX)(\ol{\exp sY \exp tX})^t).
$$
Therefore
\begin{align}
X(x_*(Y^*))=&C\left.\ppp{}{t}{s}\right|_{t=s=0}(\exp sY \exp tX \exp s\ol X^t \exp t\ol Y^t)\nnm\\
=&2C(XY+YX)=2C\left(YX+XY-\fr1m\tr(XY)I_{2m}\right)+\fr2mC(X,Y)I_{2m}\label{expl5-gausf}
\end{align}
implying that $x$ is locally strongly convex since $(X,Y)$ is positive definite.
Thus, at the origin $o$, the invariant affine metric $g_o$ is defined by:
$$
g_o(X,Y)=\left(\fr{4C}{\sqrt{2m}}\right)^{\fr2{n+2}}(X,Y)=-\fr2{mL_1}(X,Y),\quad X,Y\in {\mathfrak q}{\mathfrak h}_0(m).
$$
Since $g_o$ is positive definite and invariant by ${\rm Sp}(m)$, the invariant Riemannian metric $g$ induced by $g_o$ is exactly the affine metric of $x$. Once again we can prove that the real linear map
$$
Y\in {\mathfrak h}_0(2m)\mapsto XY+YX-\fr1m(XY)I_{2m}
$$
has a vanishing trace for each $X\in {\mathfrak h}_0(2m)$. With this fact we use \eqref{expl5-gausf} to find that $\xi=-L_1x$ holds identically, implying that $x$ is a hyperbolic affine sphere with affine mean curvature $L_1$.

Moreover, the involution map $\sigma:{\mathfrak s}{\mathfrak u}^*(2m,\bbc)\to {\mathfrak s}{\mathfrak u}^*(2m,\bbc)$ given by $\sigma(X)=-\bar X^t$ is isometric with respect to $g_o$, thus the invariant affine metric $g$ is symmetric, and the Fubini-Pick form $A_o$ of $x$ at the origin $o$ is (Definition \eqref{f-p})
\be\label{expl5-eq5}
A_o(X,Y,Z)=g_0\left(\left(XY+YX-\fr1m\tr(XY)I_{2m}\right),Z\right),\quad X,Y,Z\in {\mathfrak q}{\mathfrak h}_0(m)
\ee
which is invariant by the adjoint action ${\rm Sp}(m)$ and thus the ${\rm SU}^*(2m,\bbc)$-invariant $3$-form $A$ induced by $A_o$ is exactly the Fubini-Pick form of the hypersurface $x:M\to \bbr^{n+1}$. This shows that $x$ is an affine symmetric hypersphere.

\expl\label{expl6}\rm (\cite{bir-djo12}, \cite{lix13}) The standard embedding
$$x:M\equiv{\rm E}_{6(-26)}/{\rm F}_4\to\bbr^{27},$$
where ${\rm E}_{6(-26)}$ is the noncompact real group of type $\mathfrak{e}_6$ with the compact real form ${\rm F}_4$ of type ${\mathfrak f}_4$ as its maximal compact subgroup.

Let $\mathbb{O}$ be the space of octonions and $\mathfrak{J}$ be the set of $3\times 3$ Hermitian matrices with entries in $\mathbb{O}$, that is
$$
\mathfrak{J}=\{X=\lmx \xi_1&x_3&\bar x_2\\ \bar x_3&\xi_2&x_1\\
x_2&\bar x_1&\xi_3\rmx\in {\rm M}(3,\mathbb{O});\ \bar X^t=X\},
$$
where ${\rm M}(3,\mathbb{O})$ is the real vector space of all octonian square matrices of order $3$. Clearly $\mathfrak{J}$ is a $27$-dimensional real vector space and thus can be identified with $\bbr^{27}$. On $\mathfrak{J}$, the symmetric Jordan multiplication $\circ$ and the standard inner product $(\cdot,\cdot)$ on $\mathfrak{J}$ are defined as follows:
$$
X\circ Y=\fr12(XY+YX),\quad (X, Y)=\tr(X\circ Y).
$$
Furthermore, the cross product $\times$ and the determinant function $\det$ are given by
\bea
&X\times Y=\fr12(2X\circ Y-\tr(X)Y-\tr(Y)X+(\tr(X)\tr(Y)-\tr(X\circ Y))I_3)\\
&\det(X)=\fr13(X\times X, X).
\eea

The noncompact group ${\rm E}_{6(-26)}$ is defined as the set of all determinant-preserving real linear automorphism on $\mathfrak{J}$, that is
\be\label{E6(-26)} {\rm E}_{6(-26)}=\{A\in {\rm GL}_\bbr(\mathfrak{J});\ \det(AX)=\det(X),\,\forall X\in\mathfrak{J}\}.
\ee
The maximal compact subgroup of ${\rm E}_{6(-26)}$ is given by
\begin{align}
{\rm F}_4=&\{A\in {\rm E}_{6(-26)};\ A(X\circ Y)=(AX)\circ(AY),\, \forall X,Y\in\mathfrak{J}\}\label{F4-1}\\
\equiv&\{A\in {\rm E}_{6(-26)};\ A(I_3)=I_3\}.\label{F4-2}
\end{align}

For each matrix $T\in\mathfrak{J}$, there associated an element $\td T\in {\rm E}_{6(-26)}$ defined by
$$
\td T(X):=T\circ X,\quad \forall X\in \mathfrak{J}.
$$
Define
$$
\mathfrak{m}=\{\td T;\ T\in\mathfrak{J}_0\}, \text{\ where\ }
\mathfrak{J}_0=\{T\in\mathfrak{J};\ \tr T=0\},
$$
and denote by ${\mathfrak f}_4$ the Lie algebra of ${\rm F}_4$. Then by \cite{yok09}, the Lie algebra $\mathfrak{e}_{6(-26)}$ has a canonical direct decomposition as
\be\label{dec1}\mathfrak{e}_{6(-26)}=\mathfrak{f}_4+\mathfrak{m}\ee
satisfying
$[\mathfrak{f}_4,\mathfrak{m}]\subset \mathfrak{m}$, $[\mathfrak{m},\mathfrak{m}]\subset \mathfrak{f}_4$. Note that we have a natural identification $\mathfrak{m}\equiv T_oM$ where $o:={}_{I_{27}}{\rm F}_4$ with $I_{27}$ the identity element in ${\rm E}_{6(-26)}$.

Similar to the above, one can perform a computation which shows that the trace of an arbitrary element of $\mathfrak{e}_{6(-26)}$ must vanish (for the detail, see \cite{lix13}). Thus we have

{\prop\label{prop 5.1}\rm(\cite{lix13})
${\rm E}_{6(-26)}$ is a subgroup of the special linear group ${\rm SL}(27,\bbr)$.}

For any given constant $L_1<0$, set
$$C=\sqrt{3}(-3L_1)^{-14}>0$$
and then define a smooth map $f:{\rm E}_{6(-26)}\to \mathfrak{J}$ by
$f(L)=C\cdot L(I_3)$ for all $L\in {\rm E}_{6(-26)}$. Clearly, for any $L_1,L_2\in {\rm E}_{6(-26)}$, $f(L_1)=f(L_2)$ if and only if $(L_1^{-1}\circ L_2)(I_3)=I_3$. By the definition of ${\rm F}_4$, $f$ naturally induces a smooth map $x:{\rm E}_{6(-26)}/{\rm F}_4\to\bbr^{27}\equiv\mathfrak{J}$:
\be\label{e6/f4}
x({}_L{\rm F}_4)=C\cdot L(I_3),\quad \forall L\in {\rm E}_{6(-26)}.
\ee

By Proposition \ref{prop 5.1}, we can choose a volume element on $\bbr^{27}$,
say, the canonical volume element with respect to the inner product $(\cdot,\cdot)$ on $\mathfrak{J}$, so that ${\rm {\rm E}_{6(-26)}}$ can be identified with a subgroup of the
group ${\rm UA}(27)$ of unimodular affine transformation on $\bbr^{27}$.
Therefore, the induced map $x$ is equivariant as an affine hypersurface
in $\bbr^{27}$. Consequently all the equiaffine invariants of $x$ such as the affine metric, the Fubini-Pick form and the fundamental form are ${\rm {\rm E}_{6(-26)}}$-invariant.

Now for each $\td X\in\mathfrak{m}\equiv T_oM$, $X\in\mathfrak{J}_0$, $a(t):={}_{\exp t\td X}{\rm F}_4$ is a geodesic curve on $M$. It holds clearly that
$$
x_*(\td X)=\left.\dd{}{t}\right|_{t=0}x(a(t))=C\left.\dd{}{t}\right|_{t=0}(\exp t\td X(I_3))=C\td X(I_3)=C(X\circ I_3)=C\cdot X.
$$
This shows that $x$ is an immersion at $o$ and thus is an immersion globally since $x$ is equivariant. Clearly, $x$ is injective and is thus an imbedding of $M$ into $\bbr^{27}$.

Moreover, since for each $X\in \mathfrak{J}_0$,
$$(X,I_3)=\tr(X\circ I_3)=\tr X=0,$$ $x(o)$ is a transversal vector of $x$ at $o$ and thus is transversal everywhere. Furthermore,
for an arbitrary $Y\in\mathfrak{J}_0$, denote by $Y^*$ the Killing vector field on $M$ induced by $\td Y$, then the value of $Y^*$ at $a(t)$
$$
Y^*|_{a(t)}=\left.\dd{}{s}\right|_{s=0}(\,{}_{\exp s\td Y a(t)}{\rm F}_4)=\left.\dd{}{s}\right|_{s=0}({}\,_{\exp s\td Y \exp t\td X}{\rm F}_4).
$$
Therefore
$$
x_*(Y^*|_{a(t)})=C\left.\dd{}{s}\right|_{s=0}(\exp s\td Y \exp t\td X(I_3)).
$$
It follows that
\begin{align}
\td X(x_*(Y^*))=&C\left.\ppp{}{t}{s}\right|_{t=s=0}(\exp s\cdot\td Y\exp t\td X(I_3))\nnm\\
=&C(Y\circ(X\circ I_3))=C(Y\circ X)\nnm\\
=&C\left(X\circ Y-\fr13\tr(X\circ Y)I_3\right)+\fr13C(X,Y)I_3\label{gaussf}
\end{align}
implying that $x$ is locally strongly convex since $(X,Y)=\tr (X\circ Y)$ is positive definite.

Note that the inner product $(\cdot,\cdot)$ on $\mathfrak{J}_0$ is $\mathfrak{f}_4$-invariant and that the correspondence $\ \widetilde{\ }:\mathfrak{J}_0\to \mathfrak{m}$ is $\mathfrak{f}_4$-equivariant. It follows that the affine metric $g$ of $x$ is the invariant metric on ${\rm E}_{6(-26)}/{\rm F}_4$ induced by
$$
g_o(\td X,\td Y):=\left(\fr1{\sqrt{3}}C\right)^{\fr1{14}}(X,Y)=-\fr1{3L_1}(X,Y),\quad \forall\, X,Y\in\mathfrak{J}_0.
$$
Clearly, $g$ is symmetric since $g_o$ is invariant by the involution $\sigma(\td X)=-\widetilde{\ol X^t}$, $X\in {\mathfrak J}_0$.

A direct computation shows once more that for each $\td X$ with $X\in \mathfrak{J}_0$, the real linear map
$$
\td Y\mapsto \left(X\circ Y-\fr13\tr(X\circ Y)I_3\right)^{\widetilde{}},\quad\forall\ Y\in {\mathfrak J}_0
$$
is traceless. Taking the trace of \eqref{gaussf} respect to the metric $g$ and using $\hat\nabla_{\td X}\td Y^*=0$, $X,Y\in {\mathfrak J}_0$, with $\hat\nabla$ the Levi-Civita connection of $g$ and $\td Y^*$ the Killing vector field induced by $\td Y$, we find that the affine normal $\xi=-L_1\cdot x$ at $o$ and thus at everywhere. It follows that $x$ is a hyperbolic affine hypersphere with the affine mean curvature being the given number $L_1$.

On the other hand, the invariant Fubini-Pick form $A$ of $x$ is induced by the following $\mathfrak{f}_4$-invariant form $A_o$ (see \eqref{f-p} and \eqref{gaussf})
$$
A_o(\td X,\td Y,\td Z)=g_o\left(\left(X\circ Y-\fr13\tr(X\circ Y)I_3\right)^{\widetilde{}},\td Z\right),\,
\forall X,Y,Z\in \mathfrak{J}_0,
$$
where once again we have used the fact that $\hat\nabla_{\td X}Y^*=0$. In particular, $x$ is an affine symmetric hypersphere in $\bbr^{27}$.

\section{Proof of the main theorem with an application}

In this section we are going to prove the main theorem of this paper. After this we shall prove a proposition which makes it clear that our classification is essentially equivalent to a previous important one given by Z.J. Hu, H.Z. Li and L. Vrancken in \cite{hu-li-vra11}. Thus in a sense we in fact provide a direct way with shorter argument of proving the complete classification of the locally strongly convex hypersurfaces with parallel Fubini-Pick form. The main idea here has been used by H. Naitoh in \cite{nai81} to classify the irreducible totally real parallel submanifolds in the projective space.

Let $x:M^n\to\bbr^{n+1}$ be a locally strongly convex hypersphere with affine metric $g$ and Fubini-Pick form $A$, and suppose that $x$ is locally affine symmetric. Then by Definition \ref{dfn afsym}, $(M^n,g)$ is locally isometric to a simply connected symmetric space $G/K$ which is necessarily complete. Without loss of generality, we can put $M^n=G/K$. Furthermore, the Fubini-Pick form $A$ of $x$ must be an element of the set $\mathcal{S}_{(M^n,g)}(L_1)$ defined in Section 2.

Denote by $\mathfrak{g}$, $\mathfrak{k}$, respectively, the Lie algebras of $G$ and $K$, and $\mathfrak{g}=\mathfrak{k}+\mathfrak{m}$ the canonical decomposition of the symmetric Lie algebra pair $(\mathfrak{g},\mathfrak{k})$. Denote by $A_o$ the value of an element $A$ of $\mathcal{S}_{(M^n,g)}(L_1)$ at the origin point $o=\,{}_eK$, where $e$ is the unit element of the Lie group $G$. Define
\be\label{S^0(c)}
{\mathcal S}^0_{(M^n,g)}(L_1)=\{\sigma=A_o;\ A\in \mathcal S_{(M^n,g)}(L_1),\  \mathfrak{k}\cdot \sigma=0\}.
\ee
Then the Fubini-Pick form $A_o$ at the origin $o$ of an affine symmetric hypersphere $x:M^n\to\bbr^{n+1}$ with affine mean curvature $L_1$ is contained in ${\mathcal S}^0_{(M^n,g)}(L_1)$.

Since locally strongly convex affine hypersurface with vanishing Fubini-Pick form $A$ must be equiaffine equivalent to one of the quadric hypersurfaces given in Example \ref{expl1} (see Proposition \ref{expl1-prop}), we can assume that $A\neq 0$ thus, by the completeness and Theorems in \cite{li-sim-zhao93}, $x$ is a hyperbolic affine hypersphere where $(M^n,g)$ is a symmetric space of noncompact type.

If {\bf this} $(M^n,g)$ is reducible as a Riemannian manifold, then by Theorem \ref{chara} $x$ is a Calabi composition of $r$ points and $s$ irreducible hyperbolic affine hyperspheres. In what follows we consider the case that $(M^n,g)$ is irreducible.

The following Lemma is crucial in our proof of Theorem \ref{main}:

{\lem\label{sect4-lem} Let $M^n$ be a simply connected irreducible symmetric space of noncompact type and set $d_M=\dim\{\sigma\in S^3(\mathfrak{m});\ \mathfrak{k}\cdot\sigma=0\}$. Then $d_M=1$ if $M^n$ is one of the following spaces and $d_M=0$ otherwise:}
$$
{\rm SL}(m,\bbr)/{\rm SO}(m),\ m\geq 3;\quad {\rm SL}(m,\bbc)/{\rm SU}(m),\ m\geq 3;\quad {\rm SL}(m,\bbc)/{\rm SU}(m),\ m\geq 3;\quad {\rm E}_{6(-26)}/{\rm F}_4.
$$

\proof The argument in proving Lemma \ref{sect4-lem} is the same as the one used by H. Naitoh in \cite{nai81} (cf. the proof of Lemma 4.2 in \cite{nai81}). Let $\mathfrak{a}$ be a maximal abelian subspace in $\mathfrak{m}$ and $W$ the Weyl group of $M^n$ relative to $\mathfrak{a}$. Denote by $S^3(\mathfrak{m})$ and $S^3(\mathfrak{a})$ the vector space of all symmetric trilinear forms on $\mathfrak{m}$ and $\mathfrak{a}$, respectively. Then it is known that the vector subspace $\{\sigma\in S^3(\mathfrak{m});\ \mathfrak{k}\cdot\sigma=0\}$ is isomorphic to the vector subspace $\{\td\sigma\in S^3(\mathfrak{a});\ w\cdot\td\sigma=\td\sigma,\ \forall\,w\in W\}$ by the restriction to the subspace $\mathfrak{a}$. Since the Weyl group acts on $\mathfrak{a}$ irreducibly, all the $W$-invariant polynomials of degree $3$ are irreducible. Hence a basis of this vector subspace is given by all the fundamental $W$-invariant polynomials of degree $3$. The Weyl group $W$ for $M^n$ is of types $A_l$, $B_l$, $C_l$, $D_l$, $E_l$, ${\rm F}_4$, $G_2$ or $B_lC_l$ by the Araki's table (\cite{ara62}). Then by N. Bourbaki (\cite{bou68}), only the Weyl groups of type $A_l$ ($l\geq 2$) have one fundamental $W$-invariant polynomial of degree $3$ and the other Weyl groups have nothing. Thus the lemma follows easily. \endproof

By our previous assumption, the Fubini-Pick form $A$ of $x$ is non-vanishing, we have $$0<\dim \mathcal{S}^0_{(M^n,g)}(L_1)\leq d_M.$$ It follows that in our case $d_M=1$. Thus by Lemma \ref{sect4-lem}, $(M^n,g)$ can not be of constant sectional curvature.

{\prop\label{sect4-prop} Let $M^n$ be one of the symmetric spaces listed in Lemma \ref{sect4-lem} with symmetric metric $g$. If $\mathcal{S}^0_{(M^n,g)}(L_1)\neq\emptyset$, then the symmetric Riemannian metric $g$ is uniquely determined by the constant $L_1$ and $\mathcal{S}^0_{(M^n,g)}(L_1)$ contains only two elements $A_o,\phi\cdot A_o:=(\phi^{-1})^*A_o$ where $\phi$ is the symmetry of $(M^n,g)$ at the origin $o$.}

\proof Suppose $\td g$ is another symmetric Riemannian metric on $M^n$. Let $A\in\mathcal{S}^0_{(M^n,g)}(L_1)$ and $\td A\in\mathcal{S}^0_{(M^n,\td g)}(L_1)$. Then by \eqref{pre gaus_af sph1} we have
\begin{align}
R(X,Y)Z=&L_1(g(Y,Z)X-g(X,Z)Y)-[A(X),A(Y)](Z)\label{1}\\
\td R(X,Y)Z=&L_1(\td g(Y,Z)X-\td g(X,Z)Y)-[\td A(X),\td A(Y)](Z)\label{2}.
\end{align}
Since $M^n$ is irreducible, $\td g=\lambda^2 g$ for some positive constant $\lambda$ implying that $\td R(X,Y)Z=R(X,Y)Z$ for all $X,Y,Z\in TM^n$. Moreover, $A\neq 0$ and $\td A\neq 0$ because both $g$ and $\td g$ are not of constant sectional curvatures. On the other hand, since $\mathcal{S}^0_{(M^n,g)}(L_1)$, $\mathcal{S}^0_{(M^n,\td g)}(L_1)$ are both subsets of $\{\sigma\in S^3(\mathfrak{m});\ \mathfrak{k}\cdot\sigma=0\}$ and $d_M=1$, there is a nonzero number $\mu$ such that $\td A=\mu\cdot A$. Therefore \eqref{2} can be written as
\be\label{2'}
R(X,Y)Z=L_1\lambda^2(g(Y,Z)X-g(X,Z)Y)-\mu^2[A(X),A(Y)](Z).
\ee
Comparing \eqref{1} and \eqref{2'} we find
$$
(\mu^2-1)R(X,Y)Z=L_1(\mu^2-\lambda^2)(g(Y,Z)X-g(X,Z)Y).
$$
Note again that the metric $g$ is not of constant sectional curvature, hence $\mu^2=1$ and $\mu^2=\lambda^2$ since $L_1\neq 0$. It follows that $\td g=g$ (implying $\mathcal{S}^0_{(M^n,\td g)}(L_1)=\mathcal{S}^0_{(M^n, g)}(L_1)$) and $\td A=\pm A$,
it is easy to see that $\phi\cdot A=-A$.
\endproof

{\cor\label{sect4-cor} Let $M^n$ be one of the symmetric spaces listed in Lemma \ref{sect4-lem}. Then the symmetric affine hypersphere $x:M^n\to\bbr^{n+1}$ is unique up to affine equivalences.}

\proof If $x,\td x:M^n\to\bbr^{n+1}$ are two affine symmetric hyperspheres with a same affine mean curvature $L_1$ and with Fubini-Pick forms $A,\td A$, respectively. Then by Proposition \ref{sect4-prop}, the affine metrics of $x,\td x$ coincide and denoted as $g$. Therefore both $A_o,\td A_o$, the values of $A,\td A$ at $o$ respectively, are elements of $\mathcal{S}^0_{(M^n,g)}(L_1)$. So $\td A_o=\phi\cdot A_o$, or equivalently, $\td A_o=(\phi^{-1})^*A_o$. Consider the composition $\bar x:=x\circ\phi^{-1}$. Then the sufficient part of Theorem \ref{affine uniqueness} tells that the Fubini-Pick form $\bar A$ of $\bar x$ is subject to $\bar A=(\phi^{-1})^*A$. In particular, at $o$, we have $\bar A_o=(\phi^{-1})^*A_o=\td A_o$. Since $\bar A,\td A$ are invariant, $\bar A=\td A$ globally on $M^n$. Thus an application of the necessary part of Theorem \ref{affine uniqueness} shows that $\bar x$ and $\td x$ are equiaffine equivalent, implying that $\td x$ and $x$ are affine equivalent. \endproof

By summing up the foregoing discussions, we arrive at the completion of proving the main theorem (Theorem \ref{main}):

Let $x:M^n\to \bbr^{n+1}$ be a locally strongly convex and affine symmetric hypersphere with affine metric $g$ and affine mean curvature $L_1$.

(1) If the Fubini-Pick form $A$ of $x$ vanishes identically, then by Proposition \ref{expl1-prop}, $x$ must be one of the quadric hypersurfaces in Example \ref{expl1};

(2) If the Riemannian manifold $(M^n,g)$ is irreducible and $A\neq 0$, then by Lemma \ref{sect4-lem} and Corollary \ref{sect4-cor}, $x$ is affine equivalent to one of Examples \ref{expl3}--\ref{expl6} in Section 3;

(3) If $(M^n,g)$ is reducible, then by Theorem \ref{chara}, $x$ is affine equivalent to the Calabi composition of some $r$ points and $s$ hyperbolic affine hyperspheres listed in Examples \ref{expl1} and \ref{expl3}--\ref{expl6}, where $r+s\geq 2$.

Thus Theorem \ref{main} is proved.

Finally, to make an end of this article, we remark an alternate and simpler proof of a classification theorem originally proved by Hu et al in \cite{hu-li-vra11}. For doing this, we need the following results:

{\prop\label{sym paral} {\rm(\cite{lix13})}
A nondegenerate hypersurface $x:M^n\to\bbr^{n+1}$ is of parallel Fubini-Pick form $A$ if and only if $x$ is locally affine symmetric.}

\proof
First we suppose that the Fubini-Pick form $A$ of $x$ is parallel. Then by \cite{bok-nom-sim90}, $x$ must be an affine hypersphere. It then follows from \eqref{gaus_af sph} that the affine metric $g$ must be locally symmetric. Thus locally we can write $M^n=G/K$ and the canonical decomposition of the corresponding orthogonal symmetric pair $(\mathfrak{g},\mathfrak{k})$ is written as
${\mathfrak g}=\mathfrak{k}+\mathfrak{m}$ where the vector space $\mathfrak m$ is identified with $T_oM$. Here $o\in M^n$ is the base point given by $o=eK$ with $e$ the identity of $G$.
Note that, for all $X,Y_i\in {\mathfrak m}=T_oM$, $i=1,2,3$, the vector field $Y_i(t):=L_{\exp (tX)*}(Y_i)$ is the parallel translation of $Y_i$ along the geodesic $\gamma(t):=_{\exp (tX)}\!\!K$ (see, for example, \cite{hel01}).
Consequently we have
\begin{align}
&\dd{}{t}((L_{\exp (tX)}^*A)(Y_1,Y_2,Y_3))\nnm\\
=&\dd{}{t}(A_{\exp (tX)K}(L_{\exp (tX)*}(Y_1),
L_{\exp (tX)*}(Y_2),L_{\exp (tX)*}(Y_3)))\nnm\\
=&(\hat\nabla_{\gamma'(t)}A)(Y_1(t),Y_2(t),Y_3(t))=0,
\end{align}
where $\hat{\nabla}$ is the Levi-Civita connection of the metric $g$.
It follows that
\be\label{2.19-0}
A_{\exp (tX)K}(L_{\exp (tX)*}(Y_1),
L_{\exp (tX)*}(Y_2),L_{\exp (tX)*}(Y_3))
\ee
is constant with respect to the parameter $t$ and thus $A$ is $G$-invariant.

Conversely, we suppose that $M^n=G/K$ locally for some symmetric pair $(G,K)$ and that $A$ is $G$-invariant. Then for any $X,Y_i\in {\mathfrak m}=T_oM$, $i=1,2,3$, the function
\eqref{2.19-0}
is again a constant along the geodesic $\gamma(t)$. Therefore,
$$
(\hat\nabla_XA)(Y_1,Y_2,Y_3)=\left.\dd{}{t}\right|_{t=0} A_{\gamma(t)}(Y_1(t),Y_2(t),Y_3(t))=0,
$$
where we have once again used the fact that each $Y_i(t)$ is parallel along the geodesic $\gamma(t)$.
\endproof

{\prop\label{bok} {\rm(\cite{bok-nom-sim90})} A nondegenerate affine hypersurface with parallel Fubini-Pick form is necessarily an affine hypersphere.}

Now the following classification theorem comes readily from Theorem \ref{main}, Proposition \ref{sym paral} and Proposition \ref{bok}:

{\thm\label{cla thm} {\rm (cf. \cite{hu-li-vra11})}
Let $x:M^n\to \bbr^{n+1}$ ($n\geq 2$) be a locally strongly convex affine hypersurface with parallel Fubini-Pick form $A$. Then either of the following two cases holds:

$(1)$ With the affine metric $g$, the Riemannian manifold $(M^n,g)$ is irreducible and $x$ is locally equiaffine equivalent to

$(a)$ one of the three kinds of quadric affine spheres: Ellipsoid, elliptic paraboloid and hyperboloid; or

$(b)$ the standard embedding of the Riemannian symmetric space ${\rm SL}(m,\bbr)/{\rm SO}(m)$ into $\bbr^{n+1}$ with $n=\fr12m(m+1)-1$, $m\geq 3$;
or

$(c)$ the standard embedding of the Riemannian symmetric space ${\rm SL}(m,\bbc)/{\rm SU}(m)$ into $\bbr^{n+1}$ with $n=m^2-1$, $m\geq 3$; or

$(d)$ the standard embedding of the Riemannian symmetric space ${\rm SU}^*(2m)/{\rm Sp}(m)$ into $\bbr^{n+1}$ with $n=2m^2-m-1$, $m\geq 3$; or

$(e)$ the standard embedding of the Riemannian symmetric space ${\rm E}_{6(-26)}/{\rm F}_4$ into $\bbr^{27}$.

$(2)$ $(M^n,g)$ is reducible and $x$ is locally affine equivalent to the Calabi product of $r$ points and $s$ of the above irreducible hyperbolic affine spheres of lower dimensions, where $r$, $s$ are nonnegative integers and $r+s\geq 2$.}

\vskip 0.3in
\flushleft
Xingxiao Li\\
School of Mathematics and Information Sciences\\
Henan Normal University\\
XinXiang 453007, Henan\\
P.R.China\\
email: xxl@henannu.edu.cn

\flushleft
Guosong Zhao\\
School of Mathematics\\
Sichuan University\\
Chengdu 610064, Sichuan\\
P.R.China\\
email: gszhao@scu.edu.cn

\end{document}